\documentclass[11pt]{amsart}
\usepackage{amssymb,mathrsfs,graphicx,enumerate,color}
\usepackage{colortbl}
\definecolor{black}{rgb}{0.0, 0.0, 0.0}
\definecolor{red}{rgb}{1.0, 0.5, 0.5}

\title[   ]{$L^2$ decay for large perturbations of viscous shocks for multi-D Burgers equation}

\author[Kang]{Moon-Jin Kang}
\author[Oh]{HyeonSeop Oh}
\address[Moon-Jin Kang]{
\newline Department of Mathematical Sciences, \newline Korea Advanced Institute of Mathematical Sciences, Daejeon 34141, Korea}
\email{moonjinkang@kaist.ac.kr}
\address[HyeonSeop Oh]{
\newline Department of Mathematical Sciences, \newline Korea Advanced Institute of Mathematical Sciences, Daejeon 34141, Korea}
\email{ohs2509@kaist.ac.kr}

\newtheorem{theorem}{Theorem}[section]
\newtheorem{lemma}{Lemma}[section]

\newtheorem{remark}{Remark}[section]

\newcommand{\bbr}{\mathbb R}

\numberwithin{figure}{section}
\newcommand{\beq}{\begin{equation}}
\newcommand{\eeq}{\end{equation}}
\newcommand{\bsp}{\begin{split}}
\newcommand{\esp}{\end{split}}

\newcommand{\doubleintR}{{\int_{\mathbb{T}^{n-1}}\int_\bbr}}
\newcommand{\doubleint}{{\int_{\mathbb{T}^{n-1}} \int_0^1}}

\newcommand{\tiu}{{\tilde{u}}}

\newcommand{\s}{\sigma}

\def\eps{\varepsilon }

\newcommand\adots{\mathinner{\mkern2mu\raise1pt\hbox{.}
\mkern3mu\raise4pt\hbox{.}\mkern1mu\raise7pt\hbox{.}}}

\renewcommand{\div}{{\rm div}}

\newtheorem{theo}{Theorem}[section]

\newtheorem{lem}[theo]{Lemma}

\def\charf {\mbox{{\text 1}\kern-.30em {\text l}}}

\begin{document}
\bibliographystyle{acm}

\date{\today}


\thanks{\textbf{Acknowledgment.}  This work was supported by Samsung Science and Technology Foundation under Project Number SSTF-BA2102-01.
}

\begin{abstract}
We consider a planar viscous shock of moderate strength for a scalar viscous conservation law in multi-D. 
We consider a strictly convex flux, as a small perturbation of the Burgers flux, along the normal direction to the shock front. However, for the transversal directions, we do not have any restrictions on flux function.
We first show the contraction property for any large perturbations in $L^2$ of the planar viscous shock. If the initial $L^2$-perturbation is also in $L^1$, the large perturbation converges to zero in $L^2$ as time goes to infinity with $t^{-1/4}$ decay rate.  
The contraction and decay estimates hold up to dynamical shift. For the results, we do not impose any smallness conditions on the initial value.
This result extends the 1D case \cite{Kang-V-1} by the first author and Vasseur to the multi-dimensional case.
\end{abstract}

\maketitle \centerline{\date}


\section{Introduction}
\setcounter{equation}{0}
We consider a scalar viscous conservation law with a Lipschitz flux $F := (f_1,f_2, \cdots, f_{n})$:
\begin{align}
\begin{aligned} \label{VSCL}
&u_t + \div_x F(u)  = \Delta_x u, \\
& u(0,x) = u_0(x),
\end{aligned}
\end{align}
where $u=u(t,x)\in \bbr$ is unknown defined on 
 $t>0, x=(x_1,x')$ with  $x_1\in\bbr,\, x' \in \mathbb{T}^{n-1}:=\bbr^{n-1} / \mathbb{Z}^{n-1}, n\ge 2$ being $n-1$ dimensional flat torus. For simplicity, we use the notation $\Omega := \bbr \times \mathbb{T}^{n-1}$.\\
Consider a viscous planar shock wave $\tiu(\xi)=\tiu(x_1-\sigma t)$ such that as a traveling wave solution to \eqref{VSCL},  
\begin{align}
\begin{aligned} \label{Vshock}
&-\sigma \tiu' + f_1(\tiu)' = \tiu'',\\
& \tiu(\xi) \to u_\pm\, \text{ as } \xi \to \pm \infty,
\end{aligned}
\end{align}
where $\sigma$ is determined by the Rankine-Hugoniot condition 
\begin{equation*}
    \sigma = \frac{f_1(u_-)-f_1(u_+)}{u_- - u_+},
\end{equation*}
and the constants $u_\pm$ satisfy the Lax entropy condition $u_->u_+$. It is known that  \eqref{Vshock} admits a smooth profile $\tiu$ that is unique up to translation, and satisfies (by $f_1''>0$)
\begin{align*}
    \tiu' &= f_1(\tiu) - f_1(u_-) - \sigma(\tiu - u_-)\\
    &= (u_- - \tiu)\left(\frac{f_1(u_-) - f_1(u_+)}{u_- - u_+} - \frac{f_1(u_-) - f_1(\tiu)}{u_- - \tiu} \right) <0.
\end{align*}

Assume that the flux $f_1$ is strictly convex as a small $C^2$-perturbation of the Burgers flux such that
\begin{equation}\label{Bflux}
f_1(u)=a u^2+g(u), \quad a>0,\quad\mbox{where $g$ is a $C^2$-function satisfying $\|g''\|_{L^{\infty}(\bbr)} < \frac{2}{3}a$}.
\end{equation}

In this paper, we aim to show the contraction and decay in $L^2$ of any large perturbations of a viscous shock, up to a dynamical shift, for  \eqref{VSCL} with \eqref{Bflux}.
This extends the one-dimensional result \cite{Kang-V-1} by the first author and Vasseur to the multi-dimensional case. Our result is for large perturbations of a moderate shock, while the result in \cite{KVW} is for small perturbations of a large shock. 

Studying the contraction of any large perturbations in $L^2$-distance would be important for the physical system.
For the $L^2$ metric, we use the relative entropy defined by only one (quadratic) entropy, which is identical to the $L^2$ distance.  
Whereas, for the $L^1$-contraction by Kruzkhov \cite{K1}, we need a large family of entropies $\eta_k(u):=|u-k|, ~k\in\bbr$. However, since many physical systems of (viscous) conservation laws (including Euler, Navier-Stokes) have only one nontrivial entropy, using only one entropy for the contraction property is a remarkable program, as in \cite{KV21,KV-inven,KV-2shock,Wang22}.  \\
 For the decay estimate of large perturbation in $L^1$ distance, we refer to the result \cite {F-S} by Freist$\ddot{\mbox{u}}$hler and Serre, which provides a qualitative decay.  This result was improved by Kenig and Merle \cite{K-M}, with the uniform convergence. \\



Our main result is the following.
\begin{theorem}\label{mainthm}
Consider \eqref{VSCL} with \eqref{Bflux}. Given two constants $u_-$ and $\eps$ with $0<\eps < 8\pi (2a+\|g'' \|_{L^{\infty}(\bbr)})^{-1}$, let $u_+ = u_- - \eps$.
Let $\tiu$ be the viscous shock wave as in $\eqref{Vshock}$. \\
 Then, for any solution $u$ of \eqref{VSCL} with initial data $u_0$ satisfying $u_0-\tiu \in (L^1 \cap L^\infty) (\Omega)$, there exists a Lipschitz shift $X(t)$ such that 
\begin{equation}\label{cont}
\int_{\Omega} |u(t, x) - \tiu(x_1-\s t- X(t)) |^2 dx \leq \int_{\Omega} |u_0 - \tiu |^2 dx, \quad t>0,
\end{equation}
moreover, for all $t>0$,
\[
\|u(t, \cdot) - \tiu(\cdot-\s t- X(t)) \|_{L^2(\Omega)} = O(t^{-1/4}),
\]
and
\beq\label{dersh}
|X'(t)| =  O(t^{-1/4}).
\eeq
\end{theorem}

\begin{remark}
1. In the proof for the decay, we will show more precise quantitative estimates: 
\begin{equation}\label{timedecay}
 \|u^X - \tiu(\cdot - \sigma t)\|_{L^2(\Omega)} \leq \frac{2C C_0  n^{3/2} \| u_0 - \tiu \|_{L^2(\Omega)}}{(2\alpha)^{1/4}t^{1/4}  \| u_0 - \tiu \|_{L^2(\Omega)} + C C_0  n^{3/2}},
\end{equation}
where
\begin{align*}
C_0 &:= 1+\| u_0 - \tiu \|_{L^2(\Omega)}^2  +\| u_0 - \tiu \|_{L^1(\Omega)},\\
\alpha &:= \min \left\{ 2-\frac{2a+\|g'' \|_{L^{\infty}(\bbr)}}{2a-\|g'' \|_{L^{\infty}(\bbr)}}, 2\left(1-\frac{(2a+\|g'' \|_{L^{\infty}(\bbr)})^2 }{64\pi^2}\eps^2 \right) \right\},
\end{align*}
and $C>0$ is a positive constant depending on the endpoints $u_-, u_+$ and the interpolation inequality \eqref{GN}.
This gives the desired convergence of decay rate $t^{-1/4}$.
Then, we will show that \eqref{dersh} follows from the decay estimate and \eqref{defx} the definition of shift. \\
2. Since \eqref{dersh} implies $|X(t)| \le O(t^{3/4})$, that is, $|X(t)|$ grows sub-linearly, the shifted wave $\tiu(x_1-\s t- X(t))$ tends to the original shock $\tiu(x_1-\s t)$ time-asymptotically.
\end{remark}

To show the contraction, we use the $a$-contraction with shifts, especially with trivial weight $a=1$ as in \cite{Kang-V-1}. The method with non-trivial weight was used  in \cite{Kang19} for the contraction property for viscous scalar conservation law in 1D with general strictly convex flux. 

\section{Proof of contraction estimate}
\setcounter{equation}{0}
In this section, we show \eqref{cont}.\\
For simplicity, we use the change of variable $\xi=x_1-\sigma t$ and rewrite \eqref{VSCL} into
\begin{equation}\label{VSCLxi}
    u_t - \sigma u_\xi +\div F(u) = \Delta u,
\end{equation}
where the operators $\div$ and $\Delta$ are on $(\xi, x')$ variables.\\
As in \cite{Kang-V-1}, we use the relative entropy method with shifts.
For any strictly convex entropy $\eta$ of \eqref{VSCL}, we define the relative entropy functional by
\begin{equation*}
\eta(u|v) = \eta(u)-\eta(v) - \eta'(v)(u-v).
\end{equation*}
Moreover, we will use the notation $\mathcal{G}(\cdot|\cdot)$ to denote the relative functional of $\mathcal{G}$ (such as the relative entropy in the case that $\mathcal{G} = \eta$), i.e.,
\[
\mathcal{G}(u|v):= \mathcal{G}(u)-\mathcal{G}(v)-\mathcal{G}'(v)(u-v).
\]
Let $q(\cdot ; \cdot) = (q_1(\cdot ; \cdot),\cdots  q_n(\cdot;\cdot))$ be the flux of the relative entropy defined by
\begin{equation*}
q_i(u;v) = q_i(u)-q_i(v)-\eta'(v)(f_i(u)-f_i(v)), \, \, \ i = 1, 2, \cdots, n,
\end{equation*}
where $q = (q_1,\cdots q_n)$ is the entropy flux, i.e., $q_i' = \eta'f_i'$.

To prove the theorem \ref{mainthm}, we will use the Poincar\'e type inequality proved in \cite{KV21}.

\begin{lemma}\label{poincare}
\cite[Lemma 2.9]{KV21} For any $f:[0,1] \rightarrow \bbr$ satisfying $\int_0^1 z(1-z)|f'|^2\,dz<\infty$,
\begin{equation}
    \int_0^1 \left| f-\int_0^1 f \, dz \right|^2 \,dz \leq \frac{1}{2}\int_0^1 z(1-z)|f'|^2\,dz.
\end{equation}
\end{lemma}

First, we define a shift $X(t)$ as a solution to the ODE:
\begin{align}
\begin{aligned} \label{defx}
&\dot{X} = -\frac{2a+\|g'' \|_{L^{\infty}}}{2 \eps}\int_{\Omega}  (u(t,x) - \tiu(x_1-\s t -X(t)))\tiu'(x_1-\s t -X(t)) \, dx,\\
&X(0)=0.
\end{aligned}
\end{align}
Indeed, using the facts that $\|u\|_{L^\infty((0,\infty)\times\Omega)} \le \|u_0\|_{L^\infty(\Omega)}$ (by the maximum principle), and $\tiu', \tiu'' \in L^p(\Omega)$ for any $1\le p\le \infty$, we have the existence and uniqueness of Lipschitz solution  by the Cauchy-Lipschitz theorem.\\

For simplicity, we use the notation $u^X(t,\xi,x') = u(t,\xi+X(t),x')$, and $\tiu^{-X} = \tiu(\xi-X(t))$. 
We use the relative entropy method (as in \cite{Kang-V-1}) to have
\begin{align*}
&\frac{d}{dt} \int_{\Omega}  \eta(u|\tiu^{-X})\, d\xi\,dx'\\
&= \int_{\Omega}  (\eta'(u)-\eta'(\tiu^{-X}))\partial_t u - \eta''(\tiu^{-X})(u-\tiu^{-X})\partial_t(\tiu^{-X})\, d\xi\,dx'\\
&= \int_{\Omega}  (\eta'(u)-\eta'(\tiu^{-X}))(-\sigma u_\xi +\Delta u -\div F(u))\\
&\qquad -\eta''(\tiu^{-X})(u-\tiu^{-X})(-\dot{X}\partial_\xi\tiu^{-X}-\sigma \partial_\xi \tiu^{-X}+\Delta \tiu^{-X}-\div F(\tiu^{-X})) \, d\xi\,dx'\\
&=\dot{X} \int_{\Omega}  \eta''(\tiu^{-X})(u-\tiu^{-X})\partial_\xi\tiu^{-X}\,d\xi\,dx' \\
&\quad -\int_{\Omega}  \sigma\left((\eta'(u)-\eta'(\tiu^{-X}))u_\xi - \eta''(\tiu^{-X})(u-\tiu^{-X})\partial_\xi \tiu^{-X}\right)\,d\xi\,dx'\\
& -\int_{\Omega}  \div (q(u;\tiu^{-X}))+\eta''(\tiu^{-X})\partial_\xi\tiu^{-X}f_1(u|\tiu^{-X})  +\sum\limits_{i=2}^{n}\eta''(\tiu^{-X})\tiu_{x_i}^{-X}f_{i}(u|\tiu^{-X})\,d\xi\,dx'\\
&+\int_{\Omega}  (\eta'(u)-\eta'(\tiu^{-X}))\Delta u-\eta''(\tiu^{-X})(u-\tiu^{-X})\Delta \tiu^{-X}\,d\xi\,dx'.
\end{align*}
Taking $\eta(u)=u^2/2$, we get
\begin{align*}
&\frac{d}{dt} \int_{\Omega} \frac{|u-\tiu^{-X}|^2}{2}\, d\xi\,dx'\\
&=\dot{X}\int_{\Omega} (u-\tiu^{-X})\partial_\xi\tiu^{-X}\,d\xi\,dx' +\sigma \int_{\Omega} \left(\frac{|u-\tiu^{-X}|^2}{2}\right)_\xi \,d\xi\,dx'\\
&\quad -\int_{\Omega} f_1(u|\tiu^{-X})\partial_\xi\tiu^{-X}\,d\xi\,dx'+\int_{\Omega} (u-\tiu^{-X})(\Delta (u-\tiu^{-X})) \,d\xi \,dx'\\
&=\dot{X}\int_{\Omega} (u-\tiu^{-X})\partial_\xi\tiu^{-X}\,d\xi\,dx' - \int_{\Omega} f_1(u|\tiu^{-X})\partial_\xi\tiu^{-X}\,d\xi\,dx'\\
&\quad -\int_{\Omega} |\nabla (u-\tiu^{-X})|^2 \,d\xi\,dx'.
\end{align*}
Then, we use a change of variable $\xi \mapsto \xi+X(t)$ to have
\[
\frac{d}{dt} \int_{\Omega}  \eta(u|\tiu^{-X})\, d\xi\,dx' = \dot{X}Y(u^X)+B(u^X)-G(u^X),
\]
where
\begin{align*}
&Y(u^X) := \int_{\Omega} (u^X-\tiu)\tiu' \,d\xi\,dx' ,\\
&B(u^X) := \int_{\Omega}f_1(u^X|\tiu)\tiu' \,d\xi\,dx' ,\\
&G(u^X) := \int_{\Omega}  |\nabla (u^X-\tiu)|^2 \,d\xi\,dx' .
\end{align*}
Here, the notations $B, G$ represent the bad term and the good term respectively. Now, we rewrite the above functionals $Y, B, G$ and $\dot{X}$ with respect to the following variables
\begin{equation}\label{COV}
z = \frac{1}{\eps}(u_- - \tiu), \quad W = (u^X-\tiu) 
 \circ (z^{-1}, id),
\end{equation}
where $id : \mathbb{T}^{n-1} \to  \mathbb{T}^{n-1}$ is the identity operator on $ \mathbb{T}^{n-1}$.
Since $\tiu' <0$, we can use a change of variable $x \mapsto z\in[0,1]$. Moreover, note that 
\begin{equation}\label{chainrule}
\frac{dz}{d\xi} = -\frac{\tiu'}{\eps}.
\end{equation}

Since it holds from \eqref{defx} and the change of variable $x\to \xi$ that
\begin{align*}
\begin{aligned}
&\dot{X} = -\frac{2a+\|g'' \|_{L^{\infty}}}{2 \eps}\doubleintR (u(t, \xi + X(t), x') - \tiu(\xi))\tiu'(\xi) \, d\xi\,dx',
\end{aligned}
\end{align*}
we use \eqref{COV} and \eqref{chainrule} to have
\begin{equation}\label{XY}
\dot{X} = \left(\frac{2a+\|g''\|_{L^\infty}}{2}\right)\int_{\mathbb{T}^{n-1}} \overline{W} \,dx', \qquad  Y(u^X) = -\eps\int_{\mathbb{T}^{n-1}} \overline{W} \,dx',
\end{equation}
where $\overline{W} := \int_0^1 W \, dz$. \\
For the bad term $B$,
\begin{align}\label{Bad}
\begin{aligned}
B(u^X) &= -\doubleintR (a|u^X-\tiu|^2 + g(u^X|\tiu))\tiu' \,d\xi\,dx'\\
&\leq \doubleintR (a|u^X-\tiu|^2 + \frac{\|g''\|_{L^\infty}}{2}|u^X-\tiu|^2)\tiu' \,d\xi\,dx'\\
&= \left(\frac{2a+ \|g''\|_{L^\infty}}{2} \right) \eps \doubleint W^2 \,dz\,dx'.
\end{aligned}
\end{align}
Before estimating the good term $G$, note that
\begin{align*}
    \frac{dz}{d\xi} &= -\frac{1}{\eps}(f_1(\tiu)-f_1(u_-) - \sigma(\tiu - u_-))\\
    &= -\frac{1}{\eps}\left(a(\tiu^2 - u_-^2) + g(\tiu) - g(u_-) - \frac{f_1(u_+) - f_1(u_-)}{u_+ - u_-}(\tiu - u_-)\right)\\
    &=  -\frac{1}{\eps}\left(a(\tiu^2 - u_-^2) -a(\tiu - u_-)(u_+ + u_-) + g(\tiu) - g(u_-) - \frac{g(u_+) - g(u_-)}{u_+ - u_-}(\tiu - u_-)\right)\\
    &= -\frac{a}{\eps}(\tiu - u_-)(\tiu - u_+) - \frac{1}{\eps}\left(\frac{g(\tiu) - g(u_-)}{\tiu - u_-} - \frac{g(u_+) - g(u_-)}{u_+ - u_-} \right)(\tiu - u_-)\\
    &= \eps a z(1-z) -\frac{1}{\eps}\mathcal{J},
\end{align*}
where
\[
\mathcal{J}:= \left(\frac{g(\tiu) - g(u_-)}{\tiu - u_-} - \frac{g(u_+) - g(u_-)}{u_+ - u_-} \right)(\tiu - u_-) .
\]
By Taylor theorem for the function $Q(s) = \frac{g(s) - g(u_-)}{s-u_-}$, $\mathcal{J}$ can be estimated as follows:
\begin{align*}
    |\mathcal{J}| &= |Q'(u_\ast)(\tiu - u_+)(\tiu - u_-)|\\
    &= \left|-\frac{1}{2}g''(u_{\ast \ast})(\tiu -u_+)(\tiu - u_-)\right|\\
    &\leq \frac{\eps^2}{2}\|g''\|_{L^\infty}z(1-z).
\end{align*}
Then, we have
\begin{equation}\label{dzdxi}
\left(\frac{2a- \|g''\|_{L^\infty}}{2} \right)\eps z(1-z) \leq \frac{dz}{d\xi} \leq \left(\frac{2a+\|g''\|_{L^\infty} }{2} \right)\eps z(1-z).
\end{equation}
Thus, \eqref{dzdxi} yields that
\begin{align}\label{G}
\begin{aligned}
G(u^X) &= \doubleintR|\partial_x(u^X-\tiu)|^2\,d\xi\,dx' + \sum_{i=2}^{n}\doubleintR|\partial_{x_i}(u^X-\tiu)|^2\,d\xi\,dx'\\
&= \doubleint |\partial_z W|^2\frac{dz}{d\xi}\,dz\,dx' + \sum_{i=2}^{n}\doubleint |\partial_{x_i} W|^2\frac{d\xi}{dz}\,dz\,dx'\\
&\geq \left(\frac{2a - \|g''\|_{L^\infty}}{2}\right)\eps \doubleint z(1-z)|W_z|^2\,dz\,dx'  \\
&\quad + \sum_{i=2}^{n}\frac{2}{(2a+\|g''\|_{L^\infty})}\frac{1}{\eps}\doubleint \frac{|W_{x_i}|^2}{z(1-z)}\,dz\,dx'.
\end{aligned}
\end{align}
Using \eqref{XY}, \eqref{Bad}, \eqref{G}, we obtain that
\begin{align*}
&\frac{2}{(2a+\|g''\|_{L^\infty})}\frac{1}{\eps}(\dot{X}Y(u^X)+B(u^X)-G(u^X))\\ 
&\leq -\left(\int_\mathbb{T} \overline{W}\,dx'\right)^2 + \doubleint W^2 \,dz\,dx'- \frac{2a-\|g''\|_{L^\infty}}{2a+\|g''\|_{L^\infty}}\doubleint z(1-z)|W_z|^2 \,dz\,dx'\\
& \quad - \frac{4}{(2a+\|g''\|_{L^\infty})^2}\frac{1}{\eps^2}\sum_{i=2}^{n}\doubleint \frac{|W_{x_i}|^2}{z(1-z)} \, dz\,dx' \\
&= -\left(\int_{\mathbb{T}^{n-1}}\overline{W} \,dx'\right)^2 + \int_{\mathbb{T}^{n-1}} \overline{W}^2\,dx'+ \int_{\mathbb{T}^{n-1}}\underbrace{ \left(  \int_0^1W^2\,dz - \overline{W}^2  - \int_0^1 \frac{1}{2} z(1-z)|W_z|^2 \, dz \right) }_{=: J}  \,dx'\\
& \quad +\left(\frac{1}{2} - \frac{2a-\|g''\|_{L^\infty}}{2a+\|g''\|_{L^\infty}}\right)\doubleint z(1-z)|W_z|^2 \,dz\,dx' \\
&\quad - \frac{4}{(2a+\|g''\|_{L^\infty})^2}\frac{1}{\eps^2} \sum_{i=2}^{n}\doubleint \frac{|W_{x_i}|^2}{z(1-z)} \, dz\,dx' .
\end{align*}
Since it follows from Lemma \ref{poincare} that
\[
J = \int_0^1(W - \overline{W})^2  -  \frac{1}{2} \int_0^1 z(1-z)|W_z|^2 \, dz \le 0,
\]
we have
\begin{align*}
&\frac{2}{(2a+\|g''\|_{L^\infty})}\frac{1}{\eps}(\dot{X}Y(u^X)+B(u^X)-G(u^X))\\ 
&\leq \underbrace{-\left(\int_{\mathbb{T}^{n-1}}\overline{W}\,dx'\right)^2 + \int_{\mathbb{T}^{n-1}} \overline{W}^2\,dx' - \frac{4}{(2a+\|g''\|_{L^\infty})^2} \frac{1}{\eps^2} \sum_{i=2}^{n}\doubleint  \frac{|W_{x_i}|^2}{z(1-z)} \, dz\,dx' }_{=:K}\\
& \quad + \left(\frac{1}{2} - \frac{2a-\|g''\|_{L^\infty}}{2a+\|g''\|_{L^\infty}}\right)\doubleint z(1-z)|W_z|^2 \,dz\,dx'.
\end{align*}
Using $z(1-z) \leq 1/4$ on $0\leq z \leq 1$ and the Poincar\'e inequality on $\mathbb{T}^{n-1}$, we get
\begin{align*}
K &\leq \int_{\mathbb{T}^{n-1}} \left(\overline{W}- \int_{\mathbb{T}^{n-1}} \overline{W}\right)^2 \,dx'  - \frac{16}{(2a+\|g''\|_{L^\infty})^2}\frac{1}{\eps^2}\sum_{i=2}^{n} \doubleint  |W_{x_i}|^2\, dz\,dx'\\
& \leq  \frac{1}{4\pi^2} \sum_{i=2}^{n}\int_\mathbb{T}|\overline{W}_{x_i}|^2\,dx' - \frac{16}{(2a+\|g''\|_{L^\infty})^2}\frac{1}{\eps^2} \sum_{i=2}^{n} \doubleint  |W_{x_i}|^2 \, dz\,dx'\\
& \leq \left(\frac{1}{4\pi^2} -  \frac{16}{(2a+\|g''\|_{L^\infty})^2}\frac{1}{\eps^2} \right) \sum_{i=2}^{n} \doubleint |W_{x_i}|^2\,dz\,dx'. 
\end{align*}
Thus, by the assumptions $\|g''\|_{L^\infty} < \frac{2}{3}a$ and $0 < \eps < 8(2a+ \|g''\|_{L^\infty})^{-1}\pi$, we conclude that
\begin{align}\label{diffusionremain}
\begin{aligned}
&\frac{2}{(2a+\|g''\|_{L^\infty})}\frac{1}{\eps}(\dot{X}Y(u^X)+B(u^X)-G(u^X))\\
&\quad \leq \left(\frac{1}{2} - \frac{2a-\|g''\|_{L^\infty}}{2a+\|g''\|_{L^\infty}}\right)\doubleint z(1-z)|W_z|^2 \,dz\,dx'\\
&\quad \quad + \left(\frac{1}{4\pi^2} -  \frac{16}{(2a+\|g''\|_{L^\infty})^2}\frac{1}{\eps^2} \right) \sum_{i=2}^{n} \doubleint |W_{x_i}|^2\,dz\,dx'\\
&\leq 0.
\end{aligned}
\end{align}
Hence, we get the contraction \eqref{cont}.

\section{Proof of decay estimate}
\setcounter{equation}{0}

We here show \eqref{timedecay} and then \eqref{dersh}.\\
First, \eqref{diffusionremain} implies that
\begin{equation}\label{contconclusion}
\frac{d}{dt} \|u^X - \tiu\|_{L^2(\Omega)}^2 +\alpha \|\nabla(u^X - \tiu)\|_{L^2(\Omega)}^{2} \leq 0,
\end{equation}
where $\alpha>0$ is the constant defined by
\begin{equation*}
\alpha = \min \left\{ 2-\frac{2a+\|g'' \|_{L^{\infty}(\bbr)}}{2a-\|g'' \|_{L^{\infty}(\bbr)}}, 2\left(1-\frac{(2a+\|g'' \|_{L^{\infty}(\bbr)})^2 }{64\pi^2}\eps^2 \right) \right\}>0.
\end{equation*}

To get the time decay rate, we will use the Gagliardo Nirenberg type interpolation inequality in $\Omega$ which was proved by Huang and Yuan \cite{Huang21}.
\begin{lemma}
\cite[Theorem 1.4]{Huang21}
(Gagliardo Nirenberg type inquality in $\Omega := \bbr \times \mathbb{T}^{n-1}$)
Let $f \in L^1(\Omega)$ and $\nabla f \in L^2(\Omega)$, and $f$ is periodic in the $x_i$ direction for $i=2,\cdots,n$. Then, it holds that
\begin{equation}\label{GN}
\|f\|_{L^2(\Omega)} \leq \sum_{k=0}^{n-1} \| \nabla f \|_{L^2(\Omega)}^{\theta_k} \| f \|_{L^1(\Omega)}^{1-\theta_k},
\end{equation}
where $\theta_k = \frac{k+1}{k+3}$, and the constant $C>0$ is independent of $f$.
\end{lemma}

$\bullet$\textit{\bf Step 1)} To use the inequality \eqref{GN}, we need to investigate the $L^1$ bound of $u^X - \tiu$.  Note that
\begin{align*}
\|u^X - \tiu (\cdot - \sigma t)\|_{L^1(\Omega)} & \leq \|u^X - \tiu^X (\cdot - \sigma t)\|_{L^1(\Omega)} + \|\tiu^X - \tiu \|_{L^1(\Omega)}\\
& =: I_1 + I_2.
\end{align*}
For the $I_1$ term, we will use the following lemma.
\begin{lem}\label{L1cont}
\cite[Theorem 6.3.2]{dafermos2005hyperbolic} Let $u$ and $v$ be solutions to
\[
u_t +  \div_x f(u) = \Delta u,
\]
with respective initial data $u_0$ and $v_0$ that is in $L^1(\Omega) \cap L^\infty(\Omega)$. Assume that the flux $f$ is $C^1$. Then, we have the $L^1$ contraction
\[
\| u-v \|_{L^1(\Omega)} \leq \| u_0 - v_0 \|_{L^1(\Omega)}.
\]
\end{lem}

Applying the above Lemma \ref{L1cont}, we get 
\begin{align*}
I_1 &= \|u^X - \tiu^X(\cdot - \sigma t)\|_{L^1(\Omega)}\\
&= \|u - \tiu (\cdot - \sigma t)\|_{L^1(\Omega)}\\
&\leq \| u_0 - \tiu \|_{L^1(\Omega)}.
\end{align*}
Now, we will estimate the $I_2$ term. Observe that
\begin{align*}
I_2 &= \|\tiu^X - \tiu\|_{L^1(\Omega)}\\
&= \int_\bbr |\tiu(x_1 + X(t)) - \tiu(x_1)| \, dx_1.
\end{align*}
Since $\tiu$ is decreasing, $\tiu(x_1 + X(t)) - \tiu(x_1)$ has the opposite sign as $X(t)$ so that
\begin{align*}
I_2 &= -sgn(X(t))\int_\bbr (\tiu(x_1 + X(t)) - \tiu(x_1))\,d x_1 \\
&=  -sgn(X(t)) \int_\bbr \int_0^{X(t)} \partial_\zeta \tiu(x_1+\zeta)\,d\zeta\,dx_1\\
&= - \int_\bbr \int_0^{|X(t)|} \partial_\zeta \tiu(x_1+\zeta)\,d\zeta\,dx_1\\
&= - \int_0^{|X(t)|} \int_\bbr  \partial_{x_1} \tiu(x_1+\zeta) \, dx_1 \,d\zeta\\
&= |X(t)|(u_- - u_+).
\end{align*}
Thus, it suffices to estimate the $L^\infty$ bound of $X(t)$.\\
To get the $L^\infty$ bound of the shift function, consider the function
\[
F(\tau) = \int_\bbr |\tiu^\tau(x_1) - \tiu(x_1)|^2 \,dx_1.
\]
Then for any $\tau \in \bbr$,
\begin{align*}
\partial_\tau F(\tau) &=2 \int_\bbr (\tiu^\tau - \tiu) \partial_\tau(\tiu(x_1 + \tau))\,dx_1\\
&= 2 \int_\bbr (\tiu^\tau - \tiu) (\partial_{x_1}\tiu)(x_1 + \tau)\,dx_1\\
&= 2 \int_\bbr  (\partial_{x_1}\tiu)(x_1 + \tau) \int_{x_1}^{x_1 + \tau} \tiu'(\zeta)\,d\zeta \,dx_1\\
&= 2 \int_\bbr  \tiu'(x_1) \int_{x_1-\tau}^{x_1} \tiu'(\zeta)\,d\zeta \,dx_1\\
&=  2 \int_\bbr\int_{x_1-\tau}^{x_1}  \tiu'(x_1)  \tiu'(\zeta)\,d\zeta \,dx_1.
\end{align*}
We divide into three cases: $\tau > 1, \tau < -1$, and $-1 \leq \tau \leq 1$.
\begin{enumerate}
\item[(i)] For $\tau > 1$, we have
\[
\partial_\tau F(\tau) \geq 2 \int_\bbr \int_{x_1 -1}^{x_1} \tiu'(x_1)\tiu'(\zeta)\,d\zeta dx_1=:\beta_1>0.
\]
This implies that
\[
F(\tau) \geq F(1) + \beta_1(\tau - 1) \geq \beta_1(\tau - 1).
\]
Thus, we have
\[
\tau \leq \frac{F(\tau)}{\beta_1} + 1.
\]
\item[(ii)] For $\tau < -1$,
\[
\partial_\tau F(\tau) = -2 \int_\bbr \int_{x_1}^{x_1-\tau} \tiu'(x_1)\tiu'(\zeta)\,d\zeta dx_1 \leq  -2 \int_\bbr \int_{x_1}^{x_1+1}  \tiu'(x_1)\tiu'(\zeta)\,d\zeta dx_1  =: -\beta_2.
\]
So, we have
\[
-F(-\tau) \leq F(-1) - F(-\tau) \leq -\beta_2(-1-\tau).
\]
Since $F(\tau)$ is even, we get
\[
\tau \geq -\frac{F(\tau)}{\beta_2}-1.
\]
\item[(iii)] For the case $-1 \leq \tau \leq 1$, we have the trivial bound $|\tau| \leq 1$.
\end{enumerate}
Let $\beta := \min(\beta_1, \beta_2)$. Combining all three cases, we obtain that
\[
|\tau| \leq \frac{F(\tau)}{\beta}+1 = \frac{1}{\beta}\int_\bbr |\tiu^\tau- \tiu|^2\,dx_1+1.
\]
Taking $\tau \mapsto X(t)$, we have
\[
|X(t)| \leq \frac{1}{\beta}\int_\bbr |\tiu^X- \tiu|^2\,dx_1 +1 = \frac{1}{\beta}\int_\Omega (|\tiu^X(\cdot - \sigma t)- \tiu(\cdot - \sigma t)|^2+1) \,dx.
\]
Using the inequality $(a+b)^2 \geq b^2 - 2|ab|$, we finally get
\begin{align*}
|X(t)| &\leq \frac{1}{\beta}\left(\int_\Omega |(u^X - \tiu^X(x_1 - \sigma t)) + (\tiu^X(x_1 - \sigma t) - \tiu(x_1 - \sigma t))|^2 \,dx \right.\\
&\quad \, \,  \left. +\int_\Omega 2 |u^X - \tiu^X(x_1 - \sigma t)||\tiu^X(x_1 - \sigma t) - \tiu(x_1 - \sigma t)|\,dx \right) +1\\
&\leq \frac{1}{\beta}\left(\int_\Omega |(u^X - \tiu(x_1 - \sigma t))|^2\,dx + 4 \|\tiu\|_\infty \int_\Omega |u^X- \tiu^X(x_1 - \sigma t)|\,dx\right)+1\\
&\leq \frac{1}{\beta}(\| u_0 - \tiu \|_{L^2(\Omega)}^2 + 4 (u_- - u_+) \| u_0 - \tiu \|_{L^1(\Omega)})+1.
\end{align*}
Hence we get the $L^1$ estimate:
\begin{equation}\label{L1bound}
\|u^X - \tiu(\cdot - \sigma t)\|_{L^1(\Omega)} \leq C(1+\| u_0 - \tiu \|_{L^2(\Omega)}^2+\| u_0 - \tiu \|_{L^1(\Omega)}).
\end{equation}
For convenience, we put
\[
C_0 := 1+\| u_0 - \tiu \|_{L^2(\Omega)}^2  +\| u_0 - \tiu \|_{L^1(\Omega)}.
\]

$\bullet$\textit{\bf Step 2)}  By the interpolation inequality \eqref{GN}, and the $L^1$-bound \eqref{L1bound}, we have
\begin{align*}
\|u^X - \tiu(\cdot - \sigma t)\|_{L^2(\Omega)} &\leq C\sum_{k=0}^{n-1} \|\nabla(u^X - \tiu(\cdot - \sigma t))\|_{L^2(\Omega)}^{\frac{k+1}{k+3}} \|u^X - \tiu(\cdot - \sigma t)\|_{L^1(\Omega)}^{1-\frac{k+1}{k+3}}\\
& \leq C C_0 \sum_{k=0}^{n-1}(C_0^{-1}\|\nabla(u^X - \tiu(\cdot - \sigma t))\|_{L^2(\Omega)})^{\frac{k+1}{k+3}}.
\end{align*}
If we denote $A := C_0^{-1}\|\nabla(u^X - \tiu(\cdot - \sigma t))\|_{L^2(\Omega)}$, we can rewrite the above inequality as
\[
\|u^X - \tiu(\cdot - \sigma t)\|_{L^2(\Omega)} \leq CC_0 (A^{1/3} + A^{2/4} + \cdots + A^{n/(n+2)}).
\]
We claim that
\[
\|u^X - \tiu(\cdot - \sigma t)\|_{L^2(\Omega)} \leq C C_0 n A^{1/3}.
\]
To verify the claim, we divide into two cases: $A \leq 1$ and $A \geq 1$.
\begin{enumerate}
\item[(i)] Note that for $A \leq 1$, we have 
\[
\|u^X - \tiu(\cdot - \sigma t)\|_{L^2(\Omega)} \leq C C_0 n A^{1/3}.
\]
\item[(ii)] For the case $A>1$, we first have the upper bound
\[
\|u^X - \tiu(\cdot - \sigma t)\|_{L^2(\Omega)} \leq C C_0 n A^{n/(n+2)}.
\]
Thanks to the contraction estimate: 
\[
\|u^X - \tiu(\cdot - \sigma t)\|_{L^2(\Omega)} \leq \| u_0 - \tiu \|_{L^2(\Omega)} \leq \sqrt{C_0} \leq C_0,
\]
we obtain that
\begin{align*}
\|u^X - \tiu(\cdot - \sigma t)\|_{L^2(\Omega)} &= \|u^X - \tiu(\cdot - \sigma t)\|_{L^2(\Omega)}^{\frac{n+2}{3n}} \|u^X - \tiu(\cdot - \sigma t)\|_{L^2(\Omega)}^{1-\frac{n+2}{3n}} \\
&\leq (C C_0 n A^{\frac{n}{n+2}})^{\frac{n+2}{3n}} C_0^{1-\frac{n+2}{3n}} \\
& \leq C C_0 n^{\frac{n+2}{3n}} A^{1/3}\\
&\leq C C_0 n A^{1/3}.
\end{align*}
\end{enumerate}
Thus, we get
\begin{equation}\label{L2sob}
\begin{aligned}
    \|u^X - \tiu(\cdot - \sigma t)\|_{L^2(\Omega)} &\leq C C_0 n A^{1/3}\\
    &= C C_0^{2/3}n\|\nabla(u^X - \tiu(\cdot - \sigma t))\|_{L^2(\Omega)}^{\frac{1}{3}}.
\end{aligned}
\end{equation}
Then, we obtain that
\[
\|u^X - \tiu(\cdot - \sigma t)\|_{L^2(\Omega)}^6 \leq C C_0^4 n^6 \|\nabla(u^X - \tiu(\cdot - \sigma t))\|_{L^2(\Omega)}^{2}.
\]
Thus, by the contraction estimate \eqref{contconclusion}, we obtain that
\begin{align*}
\frac{d}{dt} \|u^X - \tiu(\cdot - \sigma t)\|_{L^2(\Omega)}^2 &\leq -\alpha \|\nabla(u^X - \tiu(\cdot - \sigma t))\|_{L^2(\Omega)}^{2}\\
& \leq -\frac{\alpha}{C C_0^4 n^6} \|u^X - \tiu(\cdot - \sigma t)\|_{L^2(\Omega)}^6.
\end{align*}
Let $\mathcal{F}(t) = \|u^X - \tiu(\cdot - \sigma t)\|_{L^2(\Omega)}^2 $. The above inequality implies that
\begin{align*}
\frac{1}{\mathcal{F}(t)^2} \geq \left( \frac{2\alpha}{C C_0^4 n^6}t + \frac{1}{\mathcal{F}(0)^2} \right)  = \frac{2\alpha t \| u_0 - \tiu \|_{L^2(\Omega)}^4 + C C_0^4 n^6}{C C_0^4 n^6 \| u_0 - \tiu (\cdot)\|_{L^2(\Omega)}^4}.
\end{align*}
Using the inequality $2(x+y)^{1/4} \geq x^{1/4} + y^{1/4}$, we can conclude that
\begin{align*}
\|u^X - \tiu(\cdot - \sigma t)\|_{L^2(\Omega)} = \sqrt{\mathcal{F}(t)} &\leq \left(\frac{C C_0^4  n^6 \| u_0 - \tiu \|_{L^2(\Omega)}^4}{2\alpha t \| u_0 - \tiu \|_{L^2(\Omega)}^4 + C C_0^4 n^6} \right)^\frac{1}{4}\\
& \leq \frac{2 C^{1/4} C_0  n^{3/2} \| u_0 - \tiu \|_{L^2(\Omega)}}{(2\alpha)^{1/4}t^{1/4}  \| u_0 - \tiu \|_{L^2(\Omega)} + C^{1/4} C_0 n^{3/2}}.
\end{align*}
This completes the proof of \eqref{timedecay}.

$\bullet$\textit{\bf Step 3)} Now, we will show the decay estimate for the shift function $X(t)$. By \eqref{defx} and \eqref{timedecay}, we have
\begin{align*}
|\dot{X}(t)| &\leq \frac{2a+\|g'' \|_{L^{\infty}}}{2 \eps}\|u - \tiu(\cdot - \sigma t - X(t))\|_{L^2(\Omega)} \|\tiu' \|_{L^2(\bbr)}\\
&\lesssim \frac{1}{1+t^{1/4}}.
\end{align*}
\qed
\vspace{0.5cm}
\bibliography{reference.bib}
\end{document}